\documentclass{amsart}
\usepackage{amsthm,amsmath,amssymb,amscd}
\usepackage{amscd,amssymb}
\newcommand{\thisdate}{\today}
{
   \newtheorem{theorem}{Theorem}[subsection]
   \newtheorem{proposition}[theorem]{Proposition}
   \newtheorem{lemma}[theorem]{Lemma}

   \newtheorem{corollary}[theorem]{Corollary}

}
{\theoremstyle{definition}

}
\newcommand{\spec}{\operatorname{Spec}}

\newcommand{\QQ}{{\mathbb{Q}}}
\newcommand{\LL}{{\mathbb{L}}}

\newcommand{\ZZ}{{\mathbb{Z}}}

\newcommand{\GG}{{\mathbb{G}}}

\newcommand{\bM}{{\mathbf{M}}}
\newcommand{\bQ}{{\mathbf{Q}}}
\newcommand{\bm}{{\mathbf{m}}}
\newcommand{\bk}{{\mathbf{k}}}
\newcommand{\bmu}{{\boldsymbol{\mu}}}

\newcommand{\cB}{{\mathcal B}}
\newcommand{\cC}{{\mathcal C}}
\newcommand{\cD}{{\mathcal D}}
\newcommand{\cE}{{\mathcal E}}
\newcommand{\cF}{{\mathcal F}}
\newcommand{\cG}{{\mathcal G}}
\newcommand{\cH}{{\mathcal H}}

\newcommand{\cK}{{\mathcal K}}
\newcommand{\cL}{{\mathcal L}}
\newcommand{\cM}{{\mathcal M}}
\newcommand{\cN}{{\mathcal N}}
\newcommand{\cO}{{\mathcal O}}
\newcommand{\cP}{{\mathcal P}}

\newcommand{\cS}{{\mathcal S}}
\newcommand{\cT}{{\mathcal T}}

\newcommand{\cV}{{\mathcal V}}

\newcommand{\ce}{\cE}
\newcommand{\cc}{\cC}
\newcommand{\cf}{\cF}
\newcommand{\fp}{\mathfrak{p}}
\newcommand{\co}{\cO}
\newcommand{\cl}{\cL}
\newcommand{\cp}{\cP}
\newcommand{\cb}{\cB}
\newcommand{\tensor}{\otimes}
\newcommand{\sym}{\mathbf{Sym}}
\newcommand{\Spec}{\operatorname{Spec}}

\newcommand{\Hom}{{\operatorname{Hom}}}
\newcommand{\cHom}{{{\cH}om}}
\newcommand{\Ext}{{\operatorname{Ext}}}
\newcommand{\cExt}{{{\cE}xt}}

\newcommand{\Aut}{{\operatorname{Aut}}}

\newcommand{\sh}{{\operatorname{sh}}}

\newcommand{\dar}{\downarrow}

\newcommand{\pl}{{P_{\cl}}}
\newcommand{\mgrnbar}{\overline{\cM}^{1/r}_{g,n}}
\newcommand{\bp}{\bar{\partial}}

\newcommand{\irightarrow}{\stackrel{\sim}{\rightarrow}}

\begin{document}
\title[Twisted spin curves]{Moduli of twisted spin curves}
\author[D. Abramovich]{Dan Abramovich}
\address{Department of Mathematics\\ Boston University\\ 111 Cummington
     Street\\ Boston, MA 02215\\ USA}
\email{abrmovic@math.bu.edu}
\author[T. Jarvis]{Tyler J. Jarvis}
\thanks{Research of D.A. was partially supported by NSF grants DMS-9700520 and
     DMS-0070970.}
\thanks{Research of T.J. was partially supported by NSA 
grant MDA904-99-1-0039} 
\address{Department of Mathematics\\ Brigham Young University\\ Provo, UT
84602\\ USA}
\email{jarvis@math.byu.edu}
\date{\thisdate}

\begin{abstract}
In this note we give a new, natural construction of a compactification of the
stack of smooth $r$-spin curves, which we call the stack of stable
twisted $r$-spin curves. This stack is  identified with a special
case of a stack of twisted stable maps of Abramovich and Vistoli.
Realizations in terms of admissible $\GG_\bm$-spaces and
$\QQ$-line bundles are given as well. The infinitesimal structure
of this stack is described in a relatively straightforward manner,
similar to that of usual stable curves.

We construct representable morphisms from the stacks of stable
twisted $r$-spin curves to the stacks of stable $r$-spin curves
\cite{Jarvis-new}, and show that they are isomorphisms. Many
delicate features of $r$-spin curves, including torsion free
sheaves with power maps, arise as simple by-products of twisted
spin curves. Various constructions, such as the
$\bar\partial$-operator of Seeley and Singer \cite{Seeley-Singer}
and Witten's cohomology class \cite{Witten} go through without
complications in the setting of twisted spin curves.
\end{abstract}

\maketitle

The moduli space of smooth $r$-spin curves was compactified by the
second author, using torsion free sheaves and coherent nets of
torsion free sheaves in \cite{Jarvis} and \cite{Jarvis-new}. In
order to construct a satisfactory compactification it was
necessary in those papers to study in detail the behavior of
torsion free sheaves with an $r$-power map. The infinitesimal
properties of those compactifications are quite subtle.

The purpose of this note is to give a new, natural construction of a
compactification of the stack of smooth $r$-spin curves, which we
call the stack of stable twisted $r$-spin curves, and to describe
its properties and relations to other compactifications.

In the first section we define the stack of smooth $r$-spin curves
and give an alternate construction in terms of certain principal
$\GG_\bm$-bundles.  To generalize these to stable curves we then
recall the twisted curves of \cite{Abramovich-Vistoli} and give a
natural definition of twisted stable $r$-spin curves.  We show
that over smooth curves all of these constructions are equivalent
to each other.  Moreover, we identify the stack of  twisted stable
$r$-spin curves  with a special case of a stack of twisted stable
maps of \cite{Abramovich-Vistoli}.

In the second section we describe the infinitesimal structure of
this stack.  This is  relatively straightforward and is similar to
the infinitesimal structure of usual stable curves.

In the third and fourth section we discuss more details of the
various characterizations in the non-smooth case.  In the third
section we describe certain admissible $\GG_\bm$-spaces that
correspond to twisted stable  spin curves.  In the fourth section
we describe the compactifications  of \cite{Jarvis} and
\cite{Jarvis-new} of the stack of stable $r$-spin curves, using
coherent nets of rank-one torsion-free sheaves, and we construct
representable morphisms from the stack of stable twisted $r$-spin
curves to those compactifications, which enables us to compare
them. It turns out that our stack of twisted spin curves is
isomorphic to that of coherent nets of roots introduced in
\cite{Jarvis-new}. Many of the delicate features of coherent nets
of torsion free sheaves arise as simple by-products of twisted
spin curves.

In the fifth section we discuss the appropriate generalization of
Witten's cohomology class \cite{Witten} to the setting of twisted
spin curves.  A key element of Witten's  construction is the
$\bar\partial$-operator of Seeley and Singer \cite{Seeley-Singer},
and we discuss how the $\bar\partial$-operator may be generalized
to the case of twisted spin curves.

\section{Spin curves and twisted spin curves}

We fix an integer $r>0$, and throughout this paper we will consider
only schemes over $\Spec \ZZ[1/r]$.

\subsection{Moduli of smooth $r$-spin curves}\label{smooth}
Let $g$ and $n$ be  integers such
that $2g-2+n>0$. Let $\bm=(m_1,\ldots,m_n)$ be an $n$-tuple of
integers.
A nonsingular $n$-pointed $r$-spin curve of genus $g$ and type
$\bm-\mathbf{1}$ over a scheme $S$, denoted $(C\to S, s_i, \cL,
c)$, is the data of
\begin{enumerate}
\item a smooth,
$n$-pointed curve $(C\to S, s_i: S \to C)$ of genus $g$,
\item an invertible sheaf $\cL$ on $C$, and
\item an isomorphism $c:\cL^{\otimes r} \stackrel{\sim}{\longrightarrow}
\omega_{C/S}(\sum_{i=1}^n (1-m_i)S_i)$ where $S_i$ is the image of
$s_i$.
\end{enumerate}

In simple terms, a nonsingular $r$-spin curve of type
$\bm-\mathbf{1}$ is a pointed curve with the choice of an $r$-th
root of the sheaf of logarithmic 1-forms which vanish to order
$m_i$ along the $i$th section.  We should point out that, while
the definition given here in terms of logarithmic differentials
vanishing to order $\bm$ is the most convenient to work with when
studying twisted spin curves (due to invariance of the generators
$dz/z$), the equivalent definition in terms of regular
differentials vanishing to order $\bm-\mathbf{1}$ is more common
(compare, for example, \cite{Jarvis-new,Witten}).

A {\em morphism} from one $r$-spin curve $(C'\to S', s'_i, \cL',
c')$ to another $(C\to S, s_i, \cL, c)$ consists of {\em a fiber
diagram} $(\phi, \alpha)$, i.e., an $S$-morphism $\phi:C' \to C$,
inducing an isomorphism $C'\to C\times_SS'$, $$\begin{array}{ccc}
C' & \stackrel{\phi}{\longrightarrow} & C \\ \dar & \square &
\dar\\ S' &\longrightarrow& S
\end{array}
$$
along with an isomorphism $ \alpha: \cL' \to \phi^*\cL$, such that
$\phi^*c\circ\alpha^{\otimes r} =c'$:

The category $\cM_{g,n}^{1/r,\bm- \mathbf{1}}$ of nonsingular,
$n$-pointed $r$-spin curves of genus $g$ and type
$\bm-\mathbf{1}$, with morphisms given by fiber diagrams, is a
Deligne-Mumford stack with quasi-projective coarse moduli space
$\bM_{g,n}^{1/r,\bm-\mathbf{1}}$. See \cite{Jarvis-new} for a
detailed proof. When $\bk$ is congruent to $\bk' \bmod r$, the two
stacks $\cM_{g,n}^{1/r,\bk}$ and $\cM_{g,n}^{1/r,\bk'}$are
canonically isomorphic.  We denote by $\cM^{1/r}_{g,n}$ the
disjoint union $\displaystyle \coprod_{-1 \leq k_i <r-1}
\cM^{1/r,\bk}_{g,n}.$

Note that if $(C\to S, s_i, \cL, c)$ is a nonsingular, $n$-pointed
$r$-spin curve of genus $g$ and type $\bm-\mathbf{1}$, and $S$ is
the prime spectrum of a field, then $\deg \cL = (2g-2+n-\sum
m_i)/r$, therefore $\chi(C,\cL) = 1-g+(2g-2+n-\sum m_i)/r $. We
denote this integer by $\chi_{r,\bm}$.

\subsection{Principal $\GG_\bm$-bundles and $r$-spin curves}
We now introduce another category which is equivalent to
$\cM_{g,n}^{1/r,\bk}$ when all $k_i$ are $-1$.

First consider $\cB\GG_\bm$, the classifying stack of $\GG_\bm$.
This is an Artin algebraic stack whose objects over a scheme $S$
are principal $\GG_\bm$-bundles $P\to S$, and whose morphisms are
fiber diagrams. Recall that such a bundle is always affine. Indeed
we have $$P = \Spec_{\cO_S}\left( \mathop{\bigoplus}\limits_{i\in
\ZZ}\ \cL^i \right),$$ where $\cL$ is an invertible sheaf on $S$.
The total space of $\cL$ is $\Spec_{\cO_S}(\bigoplus_{i \geq
0}\cL^i).$ Thus, as usual, the data of $P\to S$ and $\cL$ are
interchangeable.

The $r$-th power map  $\GG_\bm\to \GG_\bm$ induces a 1-morphism
$\kappa_r:\cB\GG_\bm\to\cB\GG_\bm$ which sends $P\to S$ to
$P/\bmu_r\to S$, or equivalently $\cL$ to $\cL^r$. Here $$P/\bmu_r
= \Spec_{\cO_S}\Big(\bigoplus_{i\ \in\ r\ZZ}\ \cL^i\Big).$$

Let $$(C\to S ,s_i:S\to C, \cL, c:\cL^r \to
\omega_{C/S}(\textstyle \sum S_i))$$ be an object of the open
moduli stack $\cM_{g,n}^{1/r,\bk}$ with all $k_i = -1$. Then we
have a 1-commutative diagram $$
\begin{array}{lcccr}
  &                             &  \cB\GG_\bm & &\\
  & \cL\nearrow     &               & \searrow\kappa_r&\\
C && \stackrel{\omega_{C/S}(\sum S_i)}{\longrightarrow} &&
\cB\GG_\bm,
\end{array}
$$ where the isomorphism $\kappa_r\circ \cL \to \omega_{C/S}(\sum
S_i)$ is given by $c$. Indeed, such a diagram is {\em equivalent}
to the data of an object of $\cM_{g,n}^{1/r,\bk}$ with all $k_i =
-1 $. Thus the category $\cM_{g,n}^{1/r,\bk}$ with all $k_i = -1$
can be defined as a category whose objects are diagrams as above,
with morphism given by fiber diagrams.

\subsection{Twisted curves}\label{sec:tw-curves}
To generalize the observation above for general $\bm$, and over
nodal curves, we use the notion of twisted curves of
\cite{Abramovich-Vistoli}. Recall that a twisted, $n$-pointed,
nodal curve $(\cC\to S, \cS_i^\cC)$ is a Deligne-Mumford stack
$\cC$, flat of relative dimension 1 over $S$, and $n$ closed
substacks $\cS_i^\cC\subset \cC$, such that $\cC$ is nodal,
$\cS_i^\cC \to S$ are \'{e}tale gerbes, the coarse moduli space $(C,
S_i)$ forms a proper, $n$-pointed, nodal curve over $S$, and
$\pi:\cC \to C$ is an isomorphism away from the nodes and the
markings $S_i$.

Twisted curves over a strictly Henselian local ring $R$ have the
following local structure: let $p$ be a geometric point of $\cC$.
Denote by $\cC^\sh$ the strict Henselization of $\cC$ at $p$.
\begin{itemize}
\item
 If $p$
does not lie over a node or a marking, then $\cC$ is isomorphic to
$C$ around $p$, and therefore $\cC^\sh \simeq \Spec R[x]^\sh$.
\item If $p$ lies over a marking $\cS_i$, then there is an integer $l$ and an
isomorphism $\cC^\sh \simeq [U/\bmu_l]$, where $U=\Spec R[z]^\sh$
and $\bmu_l$ acts on $z$ via the standard representation $z
\mapsto \zeta_lz$. We also have $C^\sh \simeq \Spec R[x]^\sh$,
where $x = z^l$.

At such a point, the dualizing sheaf corresponds to the
$\bmu_l$-equivariant sheaf $\omega_U$, with generator $dz$, on
which $\bmu_l$ acts via the standard character. It follows that
$dz/z = l^{-1}(dx/x) $ is an invariant generator of
$\omega_\cC(\cS_i)$.
\item  If $p$ lies over a node,  there is again an integer $l$ and an
isomorphism $\cC^\sh \simeq [U/\bmu_l]$, where  $U=\Spec
(R[z,w]/(z w-t))^\sh$ for some $t \in R$ and $\bmu_l$ acts on
$(z,w)$ via $(z,w) \mapsto (\zeta_lz, \zeta_l^aw)$ for some $a \in
(\ZZ/l\ZZ)^{\times}$. We also have $C^\sh \simeq \Spec (R[x,y]/(x
y-t^l))^\sh$, where $x = z^l$ and $y=w^l$.

At such a point, the dualizing sheaf corresponds to the
$\bmu_l$-equivariant sheaf $\omega_U$, with invariant generator
$\nu^U_*(dz/z - dw/w)=l^{-1}\nu^C_*(dx/x - dy/y)$, where $\nu$
stands for  normalization. The sheaf of K\"{a}hler differentials
corresponds to the $\bmu_l$-equivariant sheaf $\Omega^1_U$, with
generators $dz, dw$ (with $\bmu_l$-action as above) and relation
$w dz + z dw = 0$.
\end{itemize}
 This analysis implies that $\omega_{\cC/S}(\sum \cS_i^\cC)$ is
isomorphic to the pullback of $\omega_{C/S}(\sum S_i)$. For simplicity we
denote
$\omega_{\log} = \omega_{\log}^\cC= \omega_{\cC/S}(\sum \cS_i^\cC)$.

The twisted stable curve is {\em balanced} if the action of the
stabilizer $G_p\simeq \bmu_l$ at a nodal geometric point $p$ of
$\cC$ on the tangent spaces of the two branches of $\cC$ at $p$
has complementary eigenvalues, i.e., $a = l-1$.

\subsection{Twisted $r$-spin curves}
A twisted $r$-spin curve is a 1-commutative diagram
$$
\begin{array}{lcccr}
    &           &  \cB\GG_\bm   &&\\
    & \cL\nearrow   &       & \searrow\kappa_r&\\
\cC && \stackrel{\omega_{\cC/S}(\sum S_i)}{\longrightarrow} &&  \cB\GG_\bm,
\end{array}
$$
where
\begin{enumerate}
\item $(\cC, \cS_i^\cC)$ is a {\em balanced} twisted nodal curve,
\item $(C, S_i)$ is a stable pointed curve, and
\item $\cL: \cC \to  \cB\GG_\bm$ is representable.
\end{enumerate}

Morphisms of twisted $r$-spin curves are given by 1-commutative
fiber diagrams, up to 2-isomorphisms, as in \cite[Proposition
4.2.2]{Abramovich-Vistoli}. We denote this category by
$\cB_{g,n}(\GG_\bm,\omega_{\log}^{1/r})$.

\subsection{Twisted $r$-spin curves and twisted stable maps}  We now identify
$\cB_{g,n} (\GG_\bm,\omega_{\log}^{1/r})$ with a certain stack of
twisted stable maps. Let $C_{g,n} \to \overline\cM_{g,n}$ be the
universal curve. Consider the stack
$$C_{g,n}(\omega_{\log}^{1/r})=
C_{g,n}\mathbin{\mathop{\times}\limits_{\cB\GG_\bm}} \cB\GG_\bm$$
where the fibered product is taken with respect to the morphism
$$C \stackrel{\omega_{\log}}{\longrightarrow} \cB\GG_\bm$$ on the
left, and with respect to $\kappa_r:\cB\GG_\bm \to \cB\GG_\bm$ on
the right.

The stack $C_{g,n}(\omega_{\log}^{1/r})$ is the stack of $r$-th
roots of $\omega_{C_{g,n}/\overline\cM_{g,n}}(\sum S_i)$ over
$C_{g,n}$, which is an \'{e}tale gerbe over $C_{g,n}$ banded by
$\bmu_r$. In particular it is a Deligne-Mumford stack.

A twisted $r$-spin curve gives rise to a representable map $\cC
\to C_{g,n}(\omega_{\log}^{1/r})$, which is
  a balanced twisted stable map. The homology
class of the image of the coarse curve $C$ is the class $F$ of a fiber of the
universal 
curve $C_{g,n}\to \overline\cM_{g,n}$: the family of coarse curves $C \to S$ 
gives rise 
to a moduli morphism $S \to \overline\cM_{g,n}$ and we have an isomorphism $C
\simeq S \times_{\overline\cM_{g,n}}C_{g,n}$. We thus have a base-preserving
functor from $\cB_{g,n} (\GG_\bm,\omega_{\log}^{1/r})$ to the stack
$\cK^{bal}_{g,n}(C_{g,n}(\omega_{\log}^{1/r}) / 
\overline\cM_{g,n}, F)$, the stack of balanced, $n$-pointed twisted stable maps
of genus 
$g$ and class $F$ into $C_{g,n}(\omega_{\log}^{1/r})$ {\em relative to the
base stack} $\overline\cM_{g,n}$ (see \cite{Abramovich-Vistoli}, esp. Section
8.3). The image lies in the closed substack where 
the markings of $C$ line up over the markings of $C_{g,n}$ (i.e. the inverse
image of  $\prod_{i=1}^n S_i^{C_{g,n}}$ under the evaluation map to
$C_{g,n}^n$). 
It is easy to see that the resulting functor  is an
equivalence.
 We thus have:

\begin{theorem} The category $\cB_{g,n}( \GG_\bm,\omega_{\log}^{1/r})$ is
equivalent to the closed substack of
$\cK^{bal}_{g,n}(C_{g,n}(\omega_{\log}^{1/r}) /
\overline\cM_{g,n}, F)$, where the inverse image of $S_i^{C_{g,n}}$ is
$\cS_i^C$. In particular, it is a Deligne-Mumford stack
admitting a projective coarse moduli space.
\end{theorem}

\section{Infinitesimal structure of the stack of twisted spin curves}

\subsection{Obstructions}
\begin{proposition} The stack $\cB_{g,n}( \GG_\bm,\omega_{\log}^{1/r})$ is
smooth.
\end{proposition}

{\bf Proof.} The relative cotangent complex $\LL_{\kappa_r}$ of $\kappa_r:
\cB\GG_\bm \to \cB\GG_\bm$ is trivial, therefore deformations and obstructions
of a twisted spin curve are identical to those of the underlying pointed
twisted curves.

As in \cite[Lemma 5.3.3]{Abramovich-Vistoli}, obstructions of a
pointed twisted curve $(\cC\to S, {\cS}_i^{\cC})$ are the same as
obstructions of $\cC \to S$, given by $\Ext^2(\Omega^1_\cC,
\cO_\cC)$. As in \cite{Deligne-Mumford}, the local-to-global
spectral sequence for $\Ext$ involves, in degree 2, the terms
$H^2(\cC,\cHom(\Omega^1_\cC, \cO_\cC))$,
$H^1(\cC,\cExt^1(\Omega^1_\cC, \cO_\cC))$, and
$H^0(\cC,\cExt^2(\Omega^1_\cC, \cO_\cC))$. The first vanishes
because of dimension reasons. Similarly, $\cExt^1(\Omega^1_\cC,
\cO_\cC)$ is supported in dimension 0 hence the second term
vanishes. Thus only the third term remains. We claim that the
sheaf $\cExt^2(\Omega^1_\cC, \cO_\cC))$ vanishes. This follows
since locally $\Omega^1_\cC$ has a 2-term locally free resolution.
For instance, at a node where $\cC\simeq [U/\bmu_l]$ with $U=\Spec
(k[z,w]/(z w))^\sh$, we have a $\bmu_l$-equivariant exact sequence
on $U$: $$ 0 \to \cO_U\stackrel{(z,w)}{\longrightarrow} \cO_U
\oplus \cO_U \stackrel{(dw, dz)}{\longrightarrow} \Omega^1_U \to
0,$$ with appropriate $\bmu_l$-weights, giving a locally free
resolution of $\Omega^1_\cC$.
 Thus twisted curves are unobstructed. \qed

\subsection{Deformations}
To determine the dimension of the stack of twisted spin curves, it suffices to
identify the deformation space of a twisted pointed curve.

\begin{proposition}
 The tangent space of the stack of twisted
curves at a point $(\cC, \cS_i^{\cC})$, such that the coarse
pointed curve is stable, has dimension $3g-3+n-u$, where $u$ is
the number of nodes in $\cC$ which are not balanced.
\end{proposition}

{\bf Proof.} As in \cite[Lemma 5.3.2]{Abramovich-Vistoli}, we have
an exact sequence $$\Hom(\Omega^1_\cC, \cO_\cC) \to H^0(\cC,
\oplus\cN_{\cS_i^\cC}) \to Def \to \Ext^1(\Omega^1_\cC, \cO_\cC)
\to 0,$$ where $Def$ is the tangent space of the stack. The kernel
of the map on the left is the space of infinitesimal automorphisms
of $(\cC, \cS_i^\cC)$, which is trivial by the stability
assumption.

At points where $\cS_i$ is untwisted, we have that $H^0(\cC,
\cN_{\cS_i^\cC})$ has dimension 1.  But at twisted markings the
normal space $ \cN_{\cS_i^\cC}$ has no nontrivial local sections:
locally we have an isomorphism $[U/\bmu_l] \to \cC^\sh$, where
$U=\Spec (k[z] )^\sh$ with the standard action of $\bmu_l$;
therefore, $\bmu_l$ acts on a generator $\partial/\partial z$ of
$\cN_{\cS_i^\cC}$ via the nontrivial character $\zeta_l \mapsto
\zeta_l^{-1}$.

The group $\Ext^i(\Omega^1_\cC, \cO_\cC)$  is dual to
$H^{1-i}(\cC, \Omega^1_\cC\otimes \omega_\cC)$. By \cite[Lemma
2.3.4]{Abramovich-Vistoli}, this is the same as $H^{1-i}(C,
\pi_*(\Omega^1_\cC\otimes \omega_\cC))$, where $\pi:\cC\to C$ is
the natural map. Let us compare  $\pi_*(\Omega^1_\cC\otimes
\omega_\cC)$ with $\Omega^1_C\otimes \omega_C$. These are clearly
isomorphic away from the twisted markings and the twisted nodes.

First consider a twisted marking $\cS_i$ where  we have an
isomorphism $[U/\bmu_l] \to \cC^\sh$ with $U$ as above. The action
of $\bmu_l$ on $dz$ is via the standard character, therefore the
invariant quadratic differentials are generated by $z^l (dz/z)^2 =
l^{-2} x (dx/x)^2 =l^{-2} (dx)^2/x$, where $x$ is a parameter on
$C$. That is, locally near such a marking we have
$\pi_*(\Omega^1_\cC\otimes \omega_\cC) = \Omega^1_C\otimes
\omega_C (S_i)$.

Now consider a node on $\cC$ with  an isomorphism $[U/\bmu_l] \to
\cC^\sh$ with $U = \Spec (k[z, w]/(z w) )^\sh$, and the action can
be described via $(z,w) \mapsto (\zeta_l z, \zeta_l^a w)$ for some
$a\in (\ZZ/l\ZZ)^{\times}$. The sheaf $\omega_\cC$ has an
invariant generator  $\nu_*(dz/z-dw/w)$, where $\nu$ is the
normalization. The sheaf $\Omega^1_\cC$ has sections $f(z) dz +
g(w) dw + \alpha z dw$. Invariant elements of the form $f(z) dz +
g(w) dw$ are exactly $\Omega^1_C/\operatorname{torsion}$, whereas
$ z dw$ is invariant if and only if $a=-1$, i.e., the node is
balanced!

All in all we have $\pi_*(\Omega^1_\cC\otimes
\omega_\cC) =\Omega^1_C\otimes \omega_C (\sum S_j) / \cT$,
where the sum  $\sum S_j$ is taken over the {\em twisted} markings, and the
sheaf $\cT$ is the torsion subsheaf supported at {\em unbalanced} nodes. The
Euler characteristic agrees with that of $(C, S_i)$ with the exception of the
torsion sections at these unbalanced nodes. The Proposition follows. \qed

\begin{corollary} The morphism $\cB_{g,n}( \GG_\bm,\omega_{\log}^{1/r})\to
{\overline{\cM}}_{g,n}$ is flat, proper and quasi-finite.
\end{corollary}

\section{Twisted spin curves and admissible
$\GG_\bm$-bundles}\label{twisted}
 Let $(\cC\to S, \cS_i^\cC, \cL,
c)$ be a twisted stable $r$-spin curve. We have a corresponding
principal $\GG_\bm$-bundle $\pl\stackrel{p}{\to} \cC$. Recall that
the universal principal bundle $\cP \GG_\bm \to \cB\GG_\bm$ is
representable---in fact, $\cP \GG_\bm$ is the base scheme $S$.
Since $\cL$ is representable, we have that $\pl$ is representable
by an algebraic space. We denote by $\bar p : \pl \to C$ the
composition $\pi\circ p$, where again $\pi:\cc \rightarrow C$ is
the natural map from $\cc$ to its coarse moduli $C$.

As indicated before, the map $p:\pl\to \cC$ is affine. Since the
functor $\pi_*: Mod_{\cO_\cC} \to Mod_{\cO_C}$ is exact, we have
that $\bar p_*$ is exact on quasi-coherent sheaves, which means,
by Serre's criterion, that $\bar p:\pl\to C$ is affine.
Explicitly, we have $\pl=\Spec_{\cO_\cC}(\oplus_{i\in \ZZ}\cL^i)$,
and denoting $\cE_i = \pi_*\cL^i$, we have
$$\pl=\Spec_{\cO_C}(\mathop{\oplus}\limits_{i\in \ZZ}\cE_i).$$ The
$\ZZ$-grading reflects the fact that we still have a
$\GG_\bm$-action, making $\pl$ a homogeneous $\GG_\bm$-fibration.
Of course $\pl \to C$ is not necessarily a principal bundle unless
$\pi$ is an isomorphism. Let us make the structure of $\pl$ more
explicit.

Recall that $\cL^r\simeq \omega_{\log}$. This means that
$[\pl/\bmu_r] = P_{\omega_{\log}}$ is the principal bundle
associated to $\omega_{\log}$. On the other hand, the schematic
quotient $\pl/\bmu_r = \Spec_{\cO_C}(\oplus_{i\in r\ZZ}\cE_i).$
But $\cE_{r k} = \pi_*\cL^{rk} = (\pi_*\omega_{\log}^\cC)^k =
(\omega_{\log}^C)^k.$ Thus $\pl/\bmu_r \simeq P_{\omega_{\log}^C}
\to C$ is the principal bundle associated to $\omega_{\log}^C$.

We can describe $\pl$ over twisted markings and nodes as follows:
if $\cC^\sh \simeq  [U/\bmu_l],$ then the pullback of $\pl$ to $U$
is isomorphic to $U \times \GG_\bm$ where $\bmu_l$ acts on
$\GG_\bm$ via some non-trivial character. It follows that
$\pl|_{C^{\sh}} \simeq U \times \GG_\bm / \bmu_l$. At a twisted
marking this implies that $\pl$ is nonsingular, with fiber of
multiplicity $l$. Over a node, $\pl$ is a normal crossings surface
whose double curve is the reduction of the fiber. Moreover, the
fact that $\cC$ is balanced implies that $\pl|_{C^{\sh}}$ has the
form $\Spec R[z,w,\xi, \xi^{-1}]/(z w-t)^{\sh}\ /\ \bmu_l$ where
the action of $\bmu_l$ on $z$ and $w$ has characters inverse to
each other. We call such $P$ {\em balanced}.

We can now reverse this process: define a category
$\cP_{g,n}(\omega_{\log}^{1/r})$
whose objects are $(P \to C \to S, S_i, c)$, where
\begin{enumerate}
\item $(C\to S, S_i)$ is a stable $n$-pointed curve  of genus $g$,
\item $P\to S$ is a flat family of normal crossing surfaces,
\item $P\to C$ is a balanced, homogeneous $\GG_\bm$-fibration, which is
principal on the complement of the nodes and markings, and
\item $c: P/\bmu_r \to P_{\omega^C_{\log}}$ an isomorphism.
\end{enumerate}
Morphisms are defined using fiber diagrams as usual.

We have:
\begin{proposition}
The stacks $\cP_{g,n}(\omega_{\log}^{1/r})$ and
$\cB_{g,n}(\GG_\bm,\omega_{\log}^{1/r})$  are isomorphic.
\end{proposition}

{\bf Proof.} The argument above gives a functor
$\cB_{g,n}(\GG_\bm,\omega_{\log}^{1/r}) \to \cP_{g,n}(\omega_{\log}^{1/r})$. In the
reverse direction, given $(P \to C \to S, S_i, c)$, consider the stack
$\cC = [P/\GG_\bm]$. Since $P$ has finite stabilizers, this is a
Deligne-Mumford
stack. Clearly $P \to \cC$ is a principal bundle. Also the associated morphism
$\cC \to
\cB\GG_\bm$ is representable since $P$ is.   The local structure of $P$
implies that $\cC$ is nodal and balanced.

It is easy to see that the two compositions of these functors are equivalent to
the identity.
\qed

Let us study the algebra $\oplus_{i\in \ZZ}\cE_i$ and its
constituents locally in more detail. On the open set where $\pi$
is an isomorphism, we have $\cE_i=\cL^i$, which is invertible, and
$\cE_1 = \cL $ is an $r$-th root of $\omega_{\log}$.

Consider an \'{e}tale neighborhood of a geometric point over a
twisted marking $\cS_i^\cC$, where $\cC^\sh\simeq [U/\bmu_l]$ and
$U=\Spec (R[z])^\sh$. The sheaf $\cL$ corresponds to an
equivariant invertible sheaf on $U$, also denoted $\cL$, such that
$\cL^r\simeq \omega_{\log}$ is invariant. Choose a semi-invariant
generator $s$ of $\cL$, such that $s^r = dz/z$. The group $\bmu_l$
acts on $s$ via a character $\rho_b$ given by $\zeta_l\mapsto
\zeta_l^b$.  This implies in particular that $s^i$ is a generator
of $\cL^i$ with character $\rho_b^i$. The representability
assumption means that $b\in (\ZZ/l\ZZ)^{\times}$.

A monomial section $z^a s$ of $\cL$ is invariant if and only if
$a+b \equiv 0 \mod l$, therefore an invariant section has the form
$f(z^l) z^a s = f(x) z^a s$, where $a=l-b$. In particular,
$\pi_*\cL$ is an invertible sheaf on $U/\bmu_l$ generated by
$z^as$. And $(\pi_*\cL)^r$ is generated by $(z^as)^r = x^{ar/l}
(dz/z) = l^{-1} x^{ar/l} (dx/x)$.  If we let $k=ar/l$, then we
have that $\pi_*\cL$ is an $r$-th root of $\omega_{\log}(-k S)$ on
$U/\bmu_l$.

A similar analysis holds for $\cE_i$. A section $z^a s^i$ is
invariant if $a+ib \equiv 0 \mod l$. Denote by $k_i$ the least
non-negative integer such that $k_i \equiv -i b(r/l) \mod r$. Then
$\cE_i$ is an $r$-th root of $\omega_{\log}^i(-k_iS)$ on
$U/\bmu_l$.

Consider now an \'{e}tale neighborhood of a geometric point over a
twisted node. Here $\cC^\sh\simeq [U/\bmu_l]$ and $U=\Spec
(R[z,w]/(z w-t))^\sh$. Recall that, according to our notation,
$\bmu_l$ acts on $z$ via the standard character $\rho_1$ and on
$w$ via its inverse $\rho_{l-1}$. Again we choose a semi-invariant
generator $s$ of $\cL$, with character $\rho_b$,  such that $s^r =
\nu_*(dz/z - dw/w)$, where $\nu$ denotes normalization.  Again the
representability assumption implies that $b\in \ZZ
/l\ZZ^{\times}$.  A monomial section $z^i s$ is invariant if
$i+b\equiv 0 \mod l$, and a monomial section $w^j s$ is invariant
if $-j + b\equiv 0 \mod l$. Thus if $a=l-b$, the sheaf $\cE_1$ is
the torsion-free sheaf generated by $u = z^a s$ and $v = w^b s$,
where $yu=t^av$ and $xv=t^b u$. Moreover, $u^r =l^{-1} x^{a
r/l}(dx/x)$ and $v^r = l^{-1} y^{b r/l}(dy/y)$. In the notation of
Section \ref{stable}, below, (or of \cite{Jarvis-new}) the sheaf
$\ce_1$ is locally isomorphic to $E_{a,b}$.

Again, a similar  analysis holds for $\cE_i$.

\section{Comparison of twisted spin curves and quasi-spin curves}

We recall some definitions from \cite{Jarvis} and \cite{Jarvis-new}

\subsection{Moduli of  stable $r$-quasi-spin curves}
In \cite{Jarvis}, a natural compactification of
$\cM_{g,n}^{1/r,\bk}$ using torsion free sheaves is proposed. The
crucial point is that as a smooth $r$-spin curve degenerates to a
nodal curve $C_0$, the sheaf $\cL$ must often degenerate to a
torsion-free sheaf at the nodes of $C_0$.  But, as shown in
\cite{Jarvis}, only special types of degeneration occur.

An $n$-pointed stable {\em  $r$-quasi-spin} curve of genus $g$ and
type $\bm-\mathbf{1} $ over a scheme $S$ is the data of
\begin{enumerate}
\item a stable $n$-pointed curve $(C\to S, s_i: S \to C)$ of genus $g$,
\item a torsion-free sheaf $\cF$ on $C$ with $\chi(C,\cF)= \chi_{r,\bm}$,
 and
\item a morphism $\cF^{\otimes r} \rightarrow
\omega_{C/S}(\sum_{i=1}^n (1-m_i)S_i)$, which restricts to an
isomorphism on the open set where $\cF$ is locally free, and such
that the length of the cokernel at each node where $\cF$ is not
locally free is precisely $r-1$.
\end{enumerate}

The last condition on the cokernel of $c$ is comparable to the
fact for twisted spin curves that the local index near each node
divides $r$.

Morphisms of stable $r$-quasi-spin curves of genus $g$ and type
$\bm-\mathbf{1}$  are defined using fiber diagrams as above. The
category $\overline{\mathcal{Q}}_{g,n}^{1/r,\bm-\mathbf{1}}$ of
$n$-pointed, stable, $r$-quasi-spin curve of genus $g$ and type
$\bm-\mathbf{1}$ with morphisms given by fiber diagrams is a
proper Deligne-Mumford stack with projective coarse moduli space
$\overline\bQ_{g,n}^{1/r,\bm-\mathbf{1}}$. For proofs  see
\cite{Jarvis}.

\subsection{Moduli of stable $r$-spin curves}\label{stable}
The stack $\overline{\mathcal{Q}}_{g,n}^{1/r,\bk}\to \Spec
\ZZ[1/r]$ is not smooth. To correct this, a refined structure on
nodal curves is introduced in \cite{Jarvis-new}.  This structure
involves restrictions on the way the spin structures may vary in
families and it involves additional data in the form of
intermediate roots of $\omega_{log}$.

As mentioned in Section \ref{smooth}, the usual definition
\cite{Jarvis-new} of an $r$-spin curve of type $\bm$ is given in
terms of roots of $\omega$ instead of $\omega_{log}$.  The shift
from $\omega$ to $\omega_{log}$ adjusts the type from $\bm$ to
$\bm-\mathbf{1}$, but even accounting for this difference, the two
definitions are not exactly the same.  Still, as explained in
\cite[Note 2.1]{Jarvis-new}, although these two definitions are
not identical, there is a canonical isomorphism between the
stacks.  Thus we will use the notation
$\overline{\mathcal{M}}^{1/r,\bm-\mathbf{1}}_{g,n}$ to denote what
we are calling $r$-spin curves of type $\bm-\mathbf{1}$, despite
the slight difference between this definition and that given in
\cite{Jarvis-new}.

An $n$-pointed stable $r$-spin curve of genus $g$ and type
$\bm-\mathbf{1}$ over a scheme $S$ is the data of
\begin{enumerate}
\item a stable $n$-pointed curve $(C\to S, s_i: S \to C)$ of genus $g$,
\item a torsion free sheaf $\cF_d$ on $C$ for each $d|r$,
\item an isomorphism $\cF_1 \stackrel{\sim}{\longrightarrow}
\omega_{C/S}(\sum_{i=1}^n S_i)$, and
\item for each pair $(d,e)$ of positive integers with $d|e|r$, a morphism
$c_{e,d}:\cF_e^{\otimes e/d} \rightarrow \cF_d$ which restricts to
an isomorphism away from the nodes and marked points of $C$,
\end{enumerate}
such that the following holds:
\begin{enumerate}
\item For each $d|r$ let $m_i^d$ be the unique
nonnegative integer less than $d$ such that $m_i^d \equiv m_i \mod
d$. Write $\bm^d = (m_1^d,\ldots,m_n^d)$. Then $c_{d,1}:
\cF_d^{\otimes d} \to\cF_1$ factors through $\cF_1(-\sum m_i^d
S_i)\simeq \omega_{C/S}(\sum (1-m_i^d) S_i)$, making $(C, s_i,
\cF_d, c_{d,1})$  into a $d$-quasi spin curve of type
$\bm^d-\mathbf{1}$.
\item For each  $d|r$ the morphism $c_{d,d}$ is the identity,
\item For each $d|d'|d''|r$, we have $c_{d',d}\circ c_{d'',d'}^{\otimes d'/d}
= c_{d'',d}$, and
\item For each $d|e|r$, the length of the cokernel of $c_{e,d}$ at each point
where $\cF_e$ is not locally free is precisely $(e/d)-1$.
\end{enumerate}

Finally, we place three conditions on the local structure of these
sheaves and homomorphisms at nodes where the sheaf $\cf_r$ is not
locally free:

\begin{enumerate}
\item First, when $\cf_r$ is not locally free at a node $\fp$ in a fiber
$C_s$ of $C/S$ over a geometric point $s \in S$,  there must be an
isomorphism of the  strict Henselization $\cO_{C,\fp}^{\sh}$ of
the local ring of $C$ at $\fp$ to the ring $A :=
(R[x,y]/xy-t^l)^{\sh}$, where $l$ is some integer dividing $r$,
the ring $R=\cO_{S,s}^{\sh}$  is the strict Henselization of the
base at $s$, and the element $t$ is in $R$.
\item Second, for each pair
of positive integers $i$ and $j$, such that $i+j=l$ we define the
$A$-module $E_{i,j}:=<\xi_1, \xi_2|t^j \xi_1=x\xi_2, t^i
\xi_2=y\xi_1>$, and for notational convenience when $i=j=0$ we
define $E_{0,0}$ to be the free $A$-module with the unusual
presentation $E_{0,0}:=<\xi_1, \xi_2| \xi_1=\xi_2>$, even if
$t=0$. As explained in Sections ~\ref{twisted} and
~\ref{sec:tw-curves}, if $\cO_{\cC,\fp}^{\sh} \cong \Spec
(R[z,w]/(z w-t))^\sh$, then $E_{a,b}$ corresponds to the
pushforward $\pi_* L$ of the free $\cO_{\cC,\fp}^{\sh}$-module $L
= <s>$, where $\bmu_l$ acts on $s$ via the representation
$\rho_b$. Given a choice of isomorphism $\co^{\sh}_{C,\fp}
\irightarrow A$, we require that for each $d$ dividing $r$, the
strict Henselization $\cf^{\sh}_d$ at $\fp$ of each sheaf in the
spin structure must be isomorphic to $E_{i_d,j_d}$, where $i_d$
and $j_d$ are the least non-negative integers congruent to $i_r d$
and $j_r d$ modulo $l$, respectively, and either $i_r+j_r=l$ or
$i_r=j_r=0$.

Note that while the choice of $i_r$ and $j_r$ uniquely determines
$i_d$ and $j_d$, it does not uniquely determine the isomorphism
$\cf^{\sh}_{d} \irightarrow E_{i_d,j_d}$, since each $E_{i,j}$ has
a non-trivial automorphism group.

\item Finally, the homomorphisms $c_{d,e}: \cf_{d}^{\tensor d/e} \to
\cf_{e} $ must form a compatible system of \emph{power maps}
\cite[def. 2.3.1]{Jarvis-new}. Essentially, this means that if
$\cf_r$ is the pushforward $\pi_* \cL$ for some line bundle $\cL$
on $\cC$, then $\cf_d$ is just the pushforward $\pi_* \cL^{r/d}$
of the $r/d$th power of $\cL$, and $c_{e,d}$ corresponds to the
map $\left( \pi_* \cL^{r/e} \right)^{\tensor e/d} \to
\pi_*(\cL^{r/d})$ induced by adjointness.

More specifically, this means that given $\co^{\sh}_{C,\fp}
\irightarrow A$ we require that for each $d$ dividing $r$ there is
a choice of isomorphism $\cf^{\sh}_d \irightarrow E_{i_d,j_d}$,
such that, in terms of the presentations $\cf_{d}^{\sh} \cong
E_{i_d,j_d} =<\xi_1,\xi_2 | t^{j_d} \xi_1 = x \xi_2, y \xi_1 =
t^{i_d} \xi_2>,$ and $\cf_{e}^{\sh} \cong E_{i_e,j_e} = <
\zeta_1,\zeta_2 | t^{j_e} \zeta_1 = x \zeta_2, y \zeta_1 = t^{i_e}
\zeta_2>,$  the induced map $\overline{c}_{d,e}:\sym^{d/e}
(E_{i_d,j_d}) \to E_{i_e,j_e}$ can be written in terms of the
generators $\xi_1^{(d/e) -k}\xi_2^{k}$ of $\sym^{d/e}
(E_{i_d,j_d})$ as
\begin{equation}\label{eq:powermap}
\xi_1^{(d/e)-k}\xi_2^{k} \mapsto \begin{cases} x^{u-k}t^{kj_d}
\zeta_1 & \text{ if $0 \le k \le u$} \\y^{v- d/e + k}t^{(d-k) i_d}
\zeta_2 & \text{ if $u < k \le d$}
\end{cases}.\end{equation} Here we define  $u:= ((d/e)i_d - i_e)/l$
and $v:= ( (d/e)j_d - j_e)/l$.
\end{enumerate}

In \cite{Jarvis-new} it is shown that the property of being a
compatible system of power maps is independent of the choice of $t
\in R$ and of isomorphism $\co^{\sh}_{C,\fp} \irightarrow
(R[x,y]/xy-t^l)^{\sh}$.  Moreover, the integers $i_d$ and $j_d$
are well defined.  Of course, changing one of the isomorphisms
$\cf^{\sh}_d \cong E_{i_d,j_d}$ will generally make the induced
maps fail to be power maps; and indeed, given that these maps must
be power maps, the choice of isomorphism $\cf^{\sh}_r \cong
E_{i_r,j_r}$ uniquely determines all of the remaining isomorphisms
for $d$ dividing $r$.

The category $\overline\cM_{g,n}^{1/r,\bm-\mathbf{1}}$ of
$n$-pointed, stable $r$-spin curves of genus $g$ and type
$\bm-\mathbf{1}$, with morphisms given by fiber diagrams, is a
proper Deligne-Mumford stack with projective coarse moduli space
$\overline\bM_{g,n}^{1/r,\bm-\mathbf{1}}$. It is shown in
\cite{Jarvis-new} that this stack is smooth over $\Spec \ZZ[1/r]$,
for arbitrary $r$. The proof is quite subtle, as can be expected
from the complexity of the definition.

There is a natural morphism $\overline\cM_{g,n}^{1/r,\bk}\to
\overline{\mathcal{Q}}_{g,n}^{1/r,\bk}$ obtained by forgetting
$\cF_d$ for $d\neq 1,r$. Since $\cF_d$ is torsion free, any
automorphism of this sheaf is determined by its value over the
generic points of $C$, and therefore the morphism above is
representable.

\subsection{The functor $\cB_{g,n}(\GG_\bm,\omega_{\log}^{1/r})\to\overline\cM_{g,n}^{1/r}$}
Given a twisted stable $r$-spin curve $(\cC \to S, \cS_i, \cL,
c)$, for each positive integer $d$ dividing $r$, we have torsion
free sheaves $\cF_{r/d} = \pi_*\cL^d$ on $C$, and the canonical
isomorphisms $(\cL^d)^{e/d} \to \cL^e$ induce maps $c_{e,d} :
\cF_e^{\otimes e/d} \rightarrow \cF_d$. Consider the functor from
$\cB_{g,n}(\GG_\bm,\omega_{\log}^{1/r})$ to
$\overline\cM_{g,n}^{1/r,\bk}$ which associates to  a twisted
stable $r$-spin curve $(\cC \to S, \cS_i, \cL, c)$ the data $(C\to
S, S_i,\{\cf_d\},\{c_{e,d}\} )$. Based on the local analysis given
in Section \ref{twisted}, it is easy to see that this functor
lands in $\overline\cM_{g,n}^{1/r,\bk}$ for some $\bk$ determined
by the local indices of the twisted spin curve. The general
functor lands in the disjoint union
$$\overline{\cM}^{1/r}_{g,n}:=\coprod_{-1 \leq k_i \leq r-2}
\overline{\cM}^{1/r,\bk}_{g,n}.$$  By composing with the functor
$\overline\cM_{g,n}^{1/r} \to \overline{\mathcal{Q}}_{g,n}^{1/r}$
above, we also get a functor
$\cB_{g,n}(\GG_\bm,\omega_{\log}^{1/r})\to\mathcal{Q}_{g,n}^{1/r}$.

Note that given a balanced $\GG_\bm$-fibration $P \to C$, an
automorphism $a\in \Aut_C(P\to C)$ is determined by its action on
the generic fibers of $P\to C$ over generic points of components
of $C$. The action is  thus determined by the induced action on
the sheaf $\cE_1$, or $\cF_r$. It follows that both functors above
are representable as well. Since they are quasi-finite and proper
(these stacks being quasi-finite and proper over
$\overline{\cM}_{g,n}$) they are in fact finite. They are also
birational, being compactifications of $\cM_{g,n}^{1/r,\bm}$.

Recall that $\overline\cM_{g,n}^{1/r,\bm}$ is smooth, in particular normal. It
follows that
\begin{proposition}
$\cB_{g,n}(\GG_\bm,\omega_{\log}^{1/r})\to\overline\cM_{g,n}^{1/r}$
is an isomorphism.
\end{proposition}

\subsection{Explicit realization of the inverse functor}

We may see the inverse functor $\overline{\cM}^{1/r}_{g,n}
\rightarrow \mathcal{P}_{g,n}(\omega^{1/r}_{log})$ directly, as
follows. Given a stable $r$-spin curve $(C \rightarrow S, S_i,
\{\cF_d\},\{c_{d,d'}\})$, we must define a $\mathbb{Z}$-graded
$\co_C$-algebra $\cG$. To do this we first define a product on
rank-one, torsion-free sheaves that is compatible with the tensor
product and the spin structure.

Near a geometric node $\fp$ of $C$ over $s$ in $S$, the strict
Henselization $\co_{C,\fp}^{\sh}$ is isomorphic to
$A=(R[x,y]/xy-t^l)^{\sh}$ where $R$ is the strict Henselization of
the base $R=\co^{\sh}_{S,s}$.  The sheaf $\cF_r$ corresponds to an
$A$-module $E_{u,v}:=<\xi_1, \xi_2|t^v \xi_1 =x \xi_2,t^u\xi_2=y
\xi_1>$, where $u+v=l$.

For any integers $i,j,i',j'$ such that $i+j=i'+j'=l$, we will
define a product $$E_{i,j} \otimes E_{i',j'} \rightarrow
E_{\overline{i+i'},\overline{j+j'}}.$$ Here, for any integer $a$,
we denote by $\overline{a}$ the smallest non-negative integer
congruent to $a \pmod l$.  Intuitively, one may think of these
modules as pushforwards $\pi_* \cl = <z^is, w^js>$ and $\pi_*
\cl'=<z^{i'}s', w^{j'}s'>$ of invertible sheaves $\cl$ and $\cl'$
on $\cc$, and the product as the map $\pi_* \cl \otimes \pi_* \cl'
\rightarrow \pi_* (\cl \otimes \cl')$ induced by adjointness.

More precisely, we define the product as follows.  First, if
$i+i'>l$, then $j+j'<l$, $i>j',$ and $i'>j$. Furthermore,
$$\overline{i+i'}=i+i'-l=i-j'=i'-j$$ and $$\overline{j+j'}=j+j'.$$
In this case, if $E_{i,j}=<\zeta_1, \zeta_2|\tau^j
\zeta_1=x\zeta_2, \tau^i\zeta_2 =y \zeta_1>,$ $E_{i'j'}=<\xi_1,
\xi_2 |\tau^{j'}\xi_1=x\xi_2, \tau^{i'} \xi_2=y\xi_2>$, and $
E_{\overline{i+i'},\overline{j+j'}}=<\nu_2,
\nu_2|\tau^{\overline{j+j}} \nu_2 =x\nu_2, \tau^{\overline{i+i'}}
\nu_2 =y \nu_2>$, then the product $E_{i,j}\otimes
E_{i',j'}\rightarrow E_{\overline{i+i'},\overline{j+j'}}$ is
defined on these generators as $$\zeta_1 \otimes \xi_1 \mapsto x
\nu_1, \qquad \zeta_1 \otimes \xi_2 \mapsto \tau^{j'} \nu_1,
\qquad \zeta_2 \otimes \xi_1 \mapsto \tau^j \nu_1, \text{ and }
\qquad \zeta_2 \otimes \xi_2 \mapsto \nu_2.$$

Similarly, if $i+i'<l$, then $j+j'>l$ and
$j-i'=j'-i=\overline{j+j'}$, and we define the product by
$$\zeta_1 \otimes \xi_1 \mapsto \nu_1, \qquad \zeta_1 \otimes
\xi_2 \mapsto \tau^{i'} \nu_2, \qquad \zeta_2 \otimes \xi_1
\mapsto \tau^i \nu_2, \text{ and } \qquad \zeta_2 \otimes \xi_2
\mapsto y \nu_2.$$

Finally if $i+i'=j+j'=l$, then
$\overline{i+i'}=\overline{j+j'}=0,$ $ i=j=i'=j'=l/2$, $E_{0,0}
=<\sigma>$ is free, and the product is $$\zeta_1 \otimes \xi_1
\mapsto x \sigma \qquad \zeta_1 \otimes \xi_2 \mapsto \tau^i
\sigma \qquad \zeta_2 \otimes \xi_1 \mapsto \tau^i \sigma \qquad
\zeta_2 \otimes \xi_2 \mapsto y \sigma.$$

Note that if $d$ is relatively prime to $l$, then for any non-free
$E_{i,j}$ with $i+j=l$, the $d$-th power $E_{\overline{di},
\overline{d j}}$ is also not free, i.e., $d i$ and $d j$ are
non-zero.

It is easy to see that these product maps are commutative and
associative and induce the power maps of [4] in the case of
$(E_{i,j})^{\otimes d}$.  They also agree with the standard tensor
product on the locus where $E_{i,j}$ and $E_{i',j'}$ are free.  By
this we mean the following: When $x$ is invertible, we have
$E_{i,j} = <\zeta_1>$ with $\zeta_2 = \tau^i \zeta_1 /x$, and
$E_{i',j'} = <\xi_1>$ with $\xi_2 = \tau^{i'} \xi_1 /x$.  Letting
$\nu_1 = \zeta_1 \tensor \xi_1$, we have an isomorphism $E_{i,j}
\tensor E_{i', j'} \to E_{\overline{i+i'}, \overline{j+j'}}$.  It
is easy to check that this ``standard tensor product" isomorphism
is identical to the isomorphism induced by our product maps.  When
$y$ is invertible, a similar calculation shows the power maps
again agree with the standard tensor product, and, of course,
these are \emph{all} compatible when $x$ and $y$ are both
invertible.

This local product map allows us to define our $\co_C$-algebra
$\cG$ from an $r$-spin structure $(\{\cf_r\},\{c_{d,d'}\})$ on $C$
as follows.  Let $\cG_0:=\co_C$, and for each positive divisor $d$
of $r$ let $\cG_d :=\cf_{r/d}$, and
$\cG_{d+nr}:=(\omega^C_{log})^{\otimes n} \otimes \cf_{r/d}$.  For
all remaining indices $d$, if $e =\gcd(d,r)$ and $e \neq d$, then
on the open locus where $\cG_e$ is locally free, let
$\cG_d=\cG^{\otimes d/e}_e$.  At each node $\fp$ where $\cG_e$ is
not locally free, the local structure is $C^{\sh} \cong \spec
(R[x,y]/xy-t^l)^{\sh}$, with $l$ a divisor of $r$, and we have an
isomorphism $\alpha: \cG^{\sh}_1 \irightarrow E_{i,j}$, where
$i=i_r$, $j=j_r$, and $i+j=l$.  Further, $\alpha$ induces an
isomorphism $\cG^{\sh}_e \cong E_{\overline{e i},\overline{e j}}$.
We let $\cG^{\sh}_d := E_{\overline{d i},\overline{d j}}$ and use
the gluing induced by the power maps $(E_{i,j}^{\otimes d})
\rightarrow (E_{\overline{ei},\overline{e j}})^{\otimes d/e}
\rightarrow E_{\overline{di},\overline{d j}}$ and by the tensor
product $\cG_1|^{\otimes d}_{C-\fp} \cong \cG_e|^{\otimes
d/e}_{C-\fp} \cong \cG_d|_{C-\fp}$.  These are compatible,  as
discussed above. This gives us $\cG_{e}$ for every $e \ge 0$.

The dual $(E_{i,j})^{\vee}$ is canonically isomorphic to $E_{j,i}$
and thus we also may define $\cG_{-d}=\cG^{\vee}_d$ and extend the
$\co_C$-algebra structure to $$\cG:=\bigoplus_{d \in \mathbb{Z}}
\cG_d.$$

There is, however, a potential ambiguity in the gluing in the
construction of the $\cG_e$, arising from the choice of
isomorphism $\alpha:\cG^{\sh}_1 \irightarrow E_{i,j}$. The
following lemma shows that this is not a problem.

\begin{lemma}\label{lm:ambiguity}
Given any two isomorphisms $\alpha:\cG^{\sh}_1 \irightarrow
E_{i,j}$ and $\beta:\cG^{\sh}_1 \irightarrow E_{i,j}$, the two
sheaves of algebras, as constructed above, are isomorphic.
\end{lemma}

\begin{proof}  Let $\cG^{\alpha}_d$ and $\cG^{\beta}_d$ be the
two sheaves induced from the two choices of isomorphism. First, if
the automorphism $\alpha^{-1}\beta$ is given by multiplication by
an element of $A^{\times}$, then $\cG^{\alpha}_d$ and
$\cG^{\beta}_d$ are clearly isomorphic (the difference in gluing
data is a \v{C}ech coboundary). Consequently, we may assume that
for each $e$ dividing $r$ and such that $\cG_e$ is locally free at
$\fp$, i.e., $e i$ and $e j$ are congruent to $0 \pmod l$, an
isomorphism $\cG_e \irightarrow E_{0,0} \cong A $ has been fixed,
and the automorphism $\beta \alpha^{-1}$ preserves the
homomorphisms $$\gamma_e:E_{i,j} \rightarrow E_{\overline{ei},
\overline{e j}} =E_{0,0}$$ induced from $c_{r,r/e}:\cG^{\otimes
e}_1 = \ce^{\otimes e}_r \rightarrow \ce_{r/e} = \cG_e$.

It is straightforward to see \cite[Prop. 4.2.12]{Jarvis-new} that
the automorphisms of $E_{i,j}$ that preserve $\gamma_e$ are
elements of $\bmu_e \times \bmu_e$, acting on $E_{i,j}= <\xi_1,
\xi_2| x \xi_2=t^i \xi_1, y \xi_1 =t^i \xi_2>$ in the obvious way.
Moreover, if either $t^i$ or $t^j$ is not zero, or if $t^i$ and
$t^j$ are zero, but the normalization of $C$ at $\fp$ is
connected, then automorphisms of $E_{i,j}$ preserving $\gamma_e$
must, in fact, lie in the image of the diagonal $\bmu_e
\hookrightarrow \bmu_e \times \bmu_e.$ Thus they are all induced
by multiplication by elements of $A^{\times}$, and consequently
induce isomorphic sheaves: $\cG^{\alpha}_d \irightarrow
\cG^{\beta}_d$.

Therefore, the only potential ambiguity arises when $t^i=t^j=0$,
the normalization of $C$ at $\fp$ is disconnected, and $\beta
\alpha^{-1}$ corresponds to an element $(\eta,\sigma)$ of $\bmu_r
\times \bmu_r$.  In this case, we have $t^l=t^{i+j}=0$, and
$A=(R[x,y]/xy)^{\sh}$, and thus $E_{i,j}=<\xi_1,\xi_2|x \xi_2=y
\xi_1 =0>$.  Likewise, since $\overline{di}$ and $\overline{dj}$
are non-zero, we have $E_{\overline{di}, \overline{dj}}=<\zeta_1,
\zeta_2|x\zeta_2=y\zeta_1=0>$.  But now the automorphisms $\eta^d$
and $\sigma^d$ applied to $\cG_d|_{C-\fp}$ on the two connected
components of $C-\fp$, respectively, give an isomorphism
$\cG^{\alpha}_d \irightarrow \cG^{\beta}_d$, as desired.
\end{proof}

Since the sheaves $\ce_i$ in the $r$-spin structure are (flat)
rank-one, torsion-free sheaves with $\ce_1 =\cG_r=\omega_{log}$,
 the scheme $$P_{\cG} :=
\Spec_{\co_C} (\cG)$$ is an object of
$\cp_{g,n}(\omega^{1/r}_{log})$, and it is straightforward to see
that this functor $$\overline{\mathcal{M}}^{1/r}_{g,n} \rightarrow
\cp_{g,n} (\omega^{1/r}_{log})$$ is the inverse of the composite
$$\cp_{g,n} (\omega^{1/r}_{log})\rightarrow \cb_{g,n} (\GG_\bm,
\omega^{1/r}_{log})\rightarrow
\overline{\mathcal{M}}^{1/r}_{g,n}.$$

\section{A few remarks on analytic twisted $r$-spin curves, Witten's class, and the $\bar\partial$-operator}

The stack of $r$-spin curves is interesting primarily because of
the role it plays in the generalized Witten conjecture
\cite{Witten}, relating intersection theory on $ \mgrnbar$ to the
Gelfand-Dikii ($KdV_r$) hierarchy.  This conjecture depends on the
construction of a \emph{virtual Euler class} $c(\cV)$ for the
pushforward of the universal $r$-spin bundle.  That is, if $q: \cD
\rightarrow \mgrnbar$ is the universal curve over $\mgrnbar$, and
if $\cf_r$ is the universal $r$th root  bundle (from the universal
coherent net) on $\cD$, then $c(\cV)$ is a virtual Euler class for
$\cV:=Rq_* \cf_r$.

In particular, when $R^0q_* \cf_r=0$, then $c(\cV)$ is just the
Euler class of $-R^1q_* \cf_r$, but in general $c(\cV)$ is not the
Euler class of the index bundle $Rq_* \cf_r$, even when such a
class is well-defined.

We now describe how to extend Witten's construction to twisted
$r$-spin curves.

\subsection{The $\bar{\partial}$-operator}

Witten's construction of $c(\cV)$ depends on extending the
operator $\bar{\partial}: \Omega^{0,0}(\cf_r) \rightarrow
\Omega^{0,1}(\cf_r)$ to nodal curves; that is, to the boundary of
$\mgrnbar$.  This is done by Seeley and Singer in
\cite{Seeley-Singer} for $r=2$, and the construction for general
$r$ is similar.  Details of the general construction are given in
\cite{Mochizuki}.  This construction applies as well in the case
of twisted $r$-spin curves.

\begin{theorem}
For a family of twisted curves $(\cc/S,\cS^e_i)$ and for any line
bundle $\cl$ on $\cc$, there is a continuous family of Fredholm
operators $\bar{\partial}:\Omega^{0,0}(\cl) \rightarrow
\Omega^{0,1}(\cl)$ from $\cl$-valued $(0,0)$-forms to $\cl$-valued
$(0,1)$-forms. Moreover, if $\pi:\cc \rightarrow C$ is the coarse
moduli space of $\cc$, then $\bar{\partial}$ induces the
$\bar{\partial}$ operator of Seeley and Singer on the pushforward
to $C$: $$\bar{\partial}: \Omega^{0,0}(\pi_* \cl) \rightarrow
\Omega^{0,1} (\pi_* \cl).$$
\end{theorem}

\begin{proof}
Near a point $p \in \cc$ we have $\cc^{\sh}=[U/\bmu_l]$.  Seeley
and Singer's local description of $\bar{\partial}$ clearly applies
to any line bundle on $U$.  Moreover $\bp$ is
$\bmu_l$-equivariant, and thus passes to $\cc$, as well as
inducing the usual $\bp$ operator on the coarse moduli $C$.
\end{proof}

\subsection{Witten's Class}
Witten's construction of the virtual class $c(\cV)$ on $\mgrnbar$
works almost verbatim for $\cB_{g,n} (\GG_\bm ,
\omega^{1/r}_{log})$, where his sheaf $\cT$, (our $\cf_r$) on the
universal curve over $\mgrnbar$ is replaced by the universal
$r$th-root line bundle $\cl$ on the universal curve over
$\cB_{g,n}(\GG_\bm , \omega^{1/r}_{log})$.

The only other difference is the fact that $\pi_* \cl$ is an $r$th
root of $\omega_{log}(-\sum m_i p_i)$, rather than an $r$th root
of $\omega(-\sum m_ip_i)$.

In terms of our construction, that means that each $m_i$, or
rather the local index of the twisted curve $\cC$ near each
canonical section $\cS_i$, must be at least $1$ to make the Witten
construction work.  In the case that the local index near any
$\cS_i$ is zero, we define the class $c(\cV)$ to be $0$.   This is
compatible with the axioms of \cite{JKV} (and the axiom of descent
\cite{JKV-descent}), provided those axioms hold for the Witten
class when the local index at each $\cS_i$ is greater than $0$.

An algebraic construction of Witten's class has recently been
given by Polishchuk and Vaintrob in \cite{Polishchuk-Vaintrob}.
Their construction should also work for twisted spin curves, since
the sheaves involved in that construction are all pushforwards (to
the base of the family) of sheaves induced on the coarse moduli
$C$ by pushing forward tautological sheaves from the twisted curve
$\cC$ to $C$.

 \end{document}